\documentclass[preprint,authoryear,12pt]{elsarticle}
\usepackage{amssymb}
\newcommand{\R}{\mathbb{R}}
\newcommand{\N}{\mathbb{N}}

\newtheorem{thm}{Theorem}
\newtheorem{lem}{Lemma}

\newdefinition{defi}{Definition}
\newdefinition{rmk}{Remark}
\newdefinition{ex}{Example}
\newproof{pf}{Proof}
\newproof{pot1}{Proof of Theorem \ref{thm1}}
\newproof{pot2}{Proof of Theorem \ref{thm2}}

\begin{document}

\begin{center}
{\Large\bf Optimal Dividend Payout Model with Risk Sensitive Preferences}
\end{center}

\begin{center}{\bf \large Nicole B\"auerle $^{a}$, Anna Ja\'skiewicz$^{b}$ } \end{center}
\noindent$^{a}$Department of Mathematics, Karlsruhe Institute of Technology, Karlsruhe, Germany
{\footnotesize {\it email: nicole.baeuerle@kit.edu}}\\
\noindent$^{b}$Faculty of Pure and Applied Mathematics, Wroc{\l}aw University of Science and Technology,  Wroc{\l}aw, Poland
{\footnotesize {\it email: anna.jaskiewicz@pwr.edu.pl}}\\

\noindent {\bf Abstract.}
We consider a discrete-time dividend payout problem with risk sensitive shareholders. It is assumed  that
they are equipped with a risk aversion coefficient and construct their discounted payoff with the help of the exponential premium principle.
This leads to a risk adjusted discounted cash flow of dividends. Within such a framework not only
the expected value of the dividends is  taken into account  but also their variability.
Our approach is motivated by a remark in \cite{gershiu}. We deal with the finite and infinite time horizon problems
and prove that, even in this general setting, the optimal dividend policy is a band policy.
We also show that the policy improvement algorithm can be used to obtain the optimal policy and the corresponding value function.
Next, an explicit example is provided, in which the optimal policy is shown to be of  a barrier type.
Finally, we present some numerical studies and discuss the influence of the risk sensitive parameter on the optimal dividend policy.\\


\section{Introduction}\label{sec:intro}
The dividend payout model in risk theory is a classical problem that was introduced by \cite{DeFin57}. Since
then there have been various extensions. The  goal  is to find  for the free surplus process of an insurance company,
a dividend payout strategy that maximises  the expected discounted dividends until ruin.
Typical models for the surplus process are compound Poisson processes, diffusion processes, general renewal processes or discrete time processes.
The reader is referred  to \cite{albrthonh} and \cite{avanzi}, where  excellent overviews of recent results are provided.

Up to now most of the research has been done for the risk neutral perspective,
where the expected discounted dividends until ruin are considered. Obviously this criterion does not take the variability of the dividends into account.
From the shareholders' perspective or  from an economic point of view  it would  be certainly desirable to reduce the variability of the dividends.
Risk  should be incorporated in any kind of economic decision and shareholders are in general risk averse.
In \cite{gershiu} the authors propose the problem of maximising the expected {\em utility} of  discounted dividends until ruin instead.
Such a criterion is able to model risk aversion. In \cite{schach} the authors consider the dividend problem with an exponential utility in a diffusion setting.
They show under some assumptions that there is  a time dependent optimal barrier. \cite{bj15} consider a discrete time setting
and prove the optimality of a band policy for the exponential utility and partly characterise the optimal dividend policy
in a power utility setting. To the best of our knowledge these are so far the only papers dealing with risk sensitive dividend problems.

In this paper, we treat now the discrete time setting with state space $\R_+$ like in \cite{abt} and \cite{a14}.
However, we propose a new approach, where we consider risk sensitive preferences.
Namely, the risk adjusted discounted cash flow of the shareholder is now  of the form
$$ V_t = \alpha_t + M(V_{t+1}), \ \mbox{ where }\ M(V_{t+1})=-\frac\beta\gamma \ln\Big(\mathbb{E}_t e^{-\gamma V_{t+1}} \Big),$$
$\alpha_t$ is the dividend paid at time $t$, $\beta\in(0,1)$ is a discount factor, $\gamma>0$ is the risk sensitive parameter and
$V_t$ is the risk adjusted discounted cash flow of dividends from time $t$ onwards. These preferences are not time additive in the future dividends anymore
and allow to model risk aversion.
Note that we are here concerned about the variability of each dividend paid. This is in contrast to \cite{schach} and \cite{bj15},
where the utility of the total discounted dividends is considered. For the exponential utility with discount factor $1$ both approaches are equivalent.

The risk sensitive preferences considered in this paper belong to a wider  class  of recursive preferences
studied extensively  in macroeconomics and finance. They enjoys attention, because they allow to disentangle risk attitudes
from intertemporal substitution. In particular, \cite{EZ89} and \cite{W90} motivated the use of the certainty equivalent $M.$
Its concavity (see footnote in Sect. \ref{sec:mod}) would amplify risk aversion above intertemporal substitution.
Furthermore, the concavity of $M$ would cause the agent to prefer the early resolution of uncertainty (see \cite{KP78}).
The aforementioned recursive preference functional is still analytically tractable and retains the main behaviour features of the  risk neutral
case with $M$ replaced by the expectation operator.
One of the first papers on optimal control with this risk adjusted certainty equivalent  in discrete time is \cite{HS95}. It considers special LQ-problems.
In recent years there is a growing number of papers that study various model aspects with risk adjusted certainty equivalent, see \cite{A05,OS96,T00, W93}.

Our model can be viewed as a  Markov decision process with specific transition probability and payoff functions. Therefore, it is worth mentioning that
Markov decision processes with  dynamic risk maps  and discounted costs were examined by  \cite{Ru10}.  However, his results do not imply ours, since he studied
bounded cost functions and  coherent risk measures. In particular, such a risk measure must be positively homogeneous. Further,
 \cite{sso13} generalise the paper of \cite{Ru10} to unbounded gains and the risk sensitive average reward case. However, in their approach they apply the weighted norm approach, which result in rather stringent assumptions. Moreover, they do not analyse the properties of an optimal policy. This analysis, in our case, is necessary to show
 that the optimal policy has a band structure.
\cite{bra14} considered general certainty equivalents  for the accumulated discounted payoffs.
All the aforementioned  papers deal with Bellman equations and discuss existence and uniqueness of solutions as well as optimal policies.
However, their results are not helpful in our special setting.

The main contributions of our paper  is  threefold.  First we are able to give a mathematically rigorous solution technique for these risk sensitive dividend
problems over a finite and an infinite time horizon. More precisely, we formulate a Bellman equation which allows to compute
the value function over a finite time horizon. We also show that these value functions monotonically approximate the value function
of the infinite horizon problem. The infinite horizon value function is also characterised as a fixed point of an operator on a certain set of functions.
Second we prove that a  stationary optimal  policy has a band structure. Hence, even in this more complicated risk sensitive setting,
we are able to confirm the same form of optimal dividend payout strategy as in the risk neutral case
(for the risk neutral model consult, e.g., \cite{mi62}, \cite{morill}, \cite{Ger74}, \cite{borch82}).
Third we show that the policy improvement algorithm is another feasible way to compute the value function
and the optimal dividend payout policy for the infinite time horizon. Finally, we give some numerical examples that
shed some light on the optimal policy. For a risk sensitive model with left-sided exponential distribution for the increments of the risk reserve,
we show under some assumptions on the parameters that a barrier policy is optimal.
This result generalises \cite{a14}. For a  risk sensitive model with double-exponential distribution for the increments of the risk reserve,
we compute the optimal policy for time horizon three explicitly.
We can see some surprising dependence of the barrier on the risk sensitive parameter.

The paper is organised as follows. In Section \ref{sec:mod}, we introduce the model and our notation.
The finite horizon problem is then considered in Section \ref{sec:finite} and the limit to the infinite horizon is discussed in  Section  \ref{sec:limit}.
In Section \ref{sec:value}, we characterise the value function as the unique fixed point of  some operator within a certain class of functions.
Next we show in Section \ref{sec:optimal} that an optimal dividend policy in this risk sensitive setting is a band policy.
Afterwards we prove the validity of the policy improvement algorithm in this risk sensitive case.
In Section \ref{sec:IF} we consider an example with left-sided exponential distribution for the increments of the risk reserve and show
that a barrier policy is optimal. In the last section we provide two examples,
where we compute the optimal risk sensitive dividend payout over a time horizon of three and discuss the influence of the risk sensitive parameter
on an optimal policy.

\section{The  Model}\label{sec:mod}

We  consider the classical dividend payout problem with risk sensitive recursive evaluation of the dividends, which are paid at discrete times,
say $n\in \mathbb{N}:=  1, 2,\ldots.$	Assume there is an initial surplus $x_1$ and usually $x_1 = x\in \mathbb{R}_+:=[0,+\infty)$.
Let $Z_n$ be the difference between premium income and claim
size in the $n$-th time interval and assume that $Z_1, Z_2,\ldots$ are independent and identically distributed random variables with distribution
$\nu$ on $\mathbb{R}$. At the beginning of each time interval the insurer can decide upon paying a dividend.
The dividend payment at time $n$ is denoted by $a_n.$ If  the current risk reserve at time $n\in\N,$ say
$x_n,$ is non-negative, then $a_n$ has to be non-negative and less or equal
to  $x_n$.
If $x_n<0,$ then the company is ruined and no further dividend can be paid.
Hence, the set  of admissible dividends is
$\mathbb{A}(x_n) := [0,x_n],$ if $x_n\ge 0$ and $\mathbb{A}(x_n) := \{0\},$ if $x_n <0.$
The evolution of the surplus is given by the following equation $X_{n+1} :=f(x_n,a_n,Z_{n}),$ where
$$
f(x_n,a_n,Z_{n}):= \left\{ \begin{array}{cl}  x_n-a_n+Z_{n}, & \mbox{if }x_n\ge 0\\
 x_n,& \mbox{if }x_n< 0.\end{array}\right. $$
For any $n\in\N,$ by $H_n$ we denote the set of all feasible histories of the process up to time $n,$ i.e.,
$$
h_{n} := \left\{ \begin{array}{cl}  x_1, & \mbox{if }n= 1\\
(x_1,a_1,x_2,\ldots, x_{n}),& \mbox{if }n\ge 2,\end{array}\right.
$$
where $a_k\in \mathbb{A}(x_k)$ for $k\in \N.$
A {\it dividend policy} $\pi = (\pi_{n})_{n\in \N}$ is a sequence of Borel measurable decision rules $\pi_n: H_n\mapsto
\mathbb{R}_+$ such that $\pi_n(h_n)\in \mathbb{A}(x_n).$
Let $\Lambda$ be the set of all real-valued Borel measurable mappings such that $\alpha(x) \in\mathbb{A}(x)$ for every $x\in\mathbb{R}.$
A policy $\pi= (\pi_{n})_{n\in \N}$ is called Markov, if  $\pi_n(h_n) = \alpha_n(x_n)$ for some $\alpha_n\in\Lambda,$
  every $h_n\in H_n$ and  $n\in\N$.
A Markov policy is stationary, if $\alpha_{n}=\alpha$  for some  $\alpha\in\Lambda$ and  all $n\in\N$. In this case, we
write $\pi=\alpha^\infty$. The sets of all policies, all Markov policies,  all stationary policies are denoted by $\Pi,$ $\Pi^M$ and $\Pi^S,$ respectively.

Ruin occurs as soon as the surplus gets negative. The epoch $\tau$ of ruin is defined as the smallest positive integer $n$
such that $x_n <0.$ The question arises as to how the risk sensitive insurance company will choose its dividend strategy to maximise the gain
of the shareholder.
In this paper, we shall consider the risk adjusted discounted cash flow of dividends in the  finite and  infinite time horizon,
derived with the aid of the {\it entropic risk measure } also known as the {\it exponential premium principle.}

Let  $X$ be a non-negative real-valued random variable with distribution $\mu$ defined on some probability space $(\Omega,{\cal F}, \mathbb{P}).$
 The entropic risk measure $\rho$ for  $X$ is defined as follows
$$ \rho(X) = -\frac 1\gamma\ln \Big( \int_{\mathbb{R}_+} e^{-\gamma x}\mu(dx)\Big),$$
where $\gamma>0$ is a {\it risk sensitivity parameter } known also as a risk coefficient.
 Let  $Y$ be also a  non-negative random variable defined on $(\Omega,{\cal F}, \mathbb{P})$.
 The following properties of $\rho$\footnote{Note that $\rho$ is indeed concave, i.e., $\rho(\lambda X+(1-\lambda)Y)\ge \lambda\rho(X)+(1-\lambda)\rho(Y)$
  for any $\lambda\in[0,1].$} are important and frequently used in our analysis:
\begin{enumerate}
  \item[(P1)] monotonicity, i.e., if $X\le Y$ $\Rightarrow$ $\rho(X)\le \rho(Y),$
  \item[(P2)] translation invariance, i.e., $\rho(X+x) = \rho(X)+x$ for all $x\in\R$,
    \item[(P3)] the Jensen inequality, i.e., $\rho(X)  \le \mathbb{E} X$.
  \end{enumerate}
Furthermore, observe that by  the Taylor expansions for the exponential and logarithmic functions, we can approximate
$\rho(X)$  as follows
$$
     \rho(X)\approx \mathbb{E}X-\frac\gamma 2Var X,
$$
if $\gamma>0$ is sufficiently close to $0.$ Therefore, if $X$ is a random payoff,
then  the agent who evaluates his expected payoff with the aid of the entropic risk measure,
is not only concerned about the expected value $\mathbb{E}X$ of the random payoff $X$, but also
about its variance. Further comments on the entropic risk measure can be found in e.g.,  \cite{fs} and references cited therein.
Note that in the actuarial literature this quantity was known earlier as  the exponential premium principle (see \cite{Ger74}).

\begin{rmk}
Note that $\rho$ can be interpreted as a {\em certainty equivalent} $\rho(X)=u^{-1}(\mathbb{E}u(X))$ with
$u(x)=e^{-\gamma x}$. The exponential utility is the only function which leads to translation-invariance (see \cite{m07}), a property which we use throughout our proofs.
\end{rmk}

Let $Z$ be a random variable with the distribution $\nu.$ Throughout the paper we shall assume that
\begin{itemize}
\item[(A1)]  $\mathbb{E}Z^+=\int_0^\infty z\nu(dz) <+\infty$,
\item[(A2)] $\nu(-\infty,0)>0,$
\item[(A3)] $\nu$ has a density $g$ with respect to the Lebesgue measure.
\end{itemize}
Assumption (A2) allows to avoid a trivial case, when the ruin will never occur under any policy $\pi\in \Pi.$
We are now going to maximise the risk adjusted discounted cash flow of dividends over a finite time horizon.
For $N=2$ and fixed policy  $\pi=(\pi_k)_{k\in\mathbb{N}}\in\Pi$ the corresponding value which has to be maximized is
$$\pi_1(x)+\beta \rho(\pi_2(x,\pi_1(x),X_2)),$$
where $x_1=x$ and $X_2=x-\pi_1(x)+Z_1$ is the random risk reserve on the second stage.
For $N=3$ it is
$$ \pi_1(x)+\beta \rho\Big(\pi_2(x,\pi_1(x),X_2)+\beta \rho\big(\pi_3(x,\pi_1(x),X_2,\pi_2(\cdot),X_3)\big)\Big),$$
where $X_3=X_2-\pi_2(x,\pi_1(x),X_2)+Z_2$ is the random risk reserve on the third stage.
In order to formalise this, it is common to work with operators which we will introduce next.
Fix $k\in\N$ and
\begin{equation}\label{eq:bbar}
\bar{b}:=\frac{\beta \mathbb{E} Z^+}{1-\beta}.
\end{equation}
Let us define
 \begin{eqnarray*} B(H_k) &:=& \{ v: H_k\mapsto \mathbb{R}_+\; |\; v \mbox{ is Borel measurable}, v(h_k) \le x_k+\bar{b} \mbox{ for } x\ge0,\\
   && \hspace*{1cm} v(h_k) =0   \mbox{ for } x<0 \}.
 \end{eqnarray*}
Let  $\pi=(\pi_k)_{k\in\mathbb{N}}\in\Pi$  be an arbitrary policy. For any  $v_{k+1}\in B(H_{k+1})$ and  given $h_k\in H_k$  we put
$$
\rho_{\pi_k,h_k}(v_{k+1}):=-\frac 1\gamma\ln\left(
 \int_{\R} e^{-\gamma v_{k+1}(h_k,\pi_k(h_k),f(x_k,\pi_k(h_k),z))}\nu(dz)\right).$$
 Hence,
$$
\rho_{\pi_k,h_k}(v_{k+1})= -\frac 1\gamma\ln\left(
 \int_{\pi_k(h_k)-x_k}^\infty e^{-\gamma v_{k+1}(h_k,\pi_k(h_k),x_k-\pi_k(h_k)+z)}\nu(dz)+\nu(-\infty,\pi_k(h_k)-x_k)\right),
 $$
 if $x_k\ge 0$ and $\rho_{\pi_k,h_k}(v_{k+1})=0,$ if $x_k<0.$
Furthermore,  we define the operator $L_{\pi_k}$ for functions $v_{k+1}\in B(H_{k+1})$ as follows
$$
(L_{\pi_k}v_{k+1})(h_k):=\pi_k(h_k)+\beta\rho_{\pi_k,h_k}(v_{k+1}),
$$
where $\beta\in (0,1)$ is a discount factor.
By property (P1), it follows that $L_{\pi_k}$ is monotone, i.e.,
\begin{equation}
\label{L_mon}
(L_{\pi_k}v_{k+1})(h_k)\le (L_{\pi_k}\hat{v}_{k+1})(h_k)\quad \mbox{if}
\quad v_{k+1}\le \hat{v}_{k+1}, \quad  v_{k+1}, \hat{v}_{k+1},\in B(H_{k+1}).
\end{equation}
We shall write $Lv$ instead of $(Lv)$. Moreover,  by (P2)  for any constant $\hat{b}\in \R_+$ we get that
\begin{equation}
\label{L_ogr}
0\le L_{\pi_k}(v_{k+1}+\hat{b})(h_k)=  L_{\pi_k}v_{k+1}(h_k) +\beta \hat{b}
\end{equation}
for every $h_k\in H_k$ with $k\in\mathbb{N}.$
For any initial income $x_1=x\in\R_+$ and  $N\in {\mathbb N}$ we define the $N$-stage total risk adjusted discounted cash flow of dividends by
\begin{equation}\label{J_T}
J_{N}(x,\pi):=(L_{\pi_1}\circ \ldots \circ L_{\pi_{N}} ){\bf 0}(x),
\end{equation}
where  ${\bf 0}$ is a function such that  ${\bf 0}(h_k)\equiv 0$ for every  $h_k\in H_k$ and $k\in\mathbb{N}.$
Note that we show below that $(L_{\pi_k}\circ \ldots \circ L_{\pi_{N}} ){\bf 0}(x) \in B(h_k)$, so the iteration is well-defined. Clearly, if $x<0,$ then $J_{N}(x,\pi)=0.$
For instance, if $N=2$ and $x\in\R_+,$ definition (\ref{J_T}) is read as follows
\begin{eqnarray*}
J_2(x,\pi)&=&
(L_{\pi_1}\circ L_{\pi_2}){\bf 0}(x)=L_{\pi_1}( L_{\pi_2}{\bf 0})(x)\\
&=& \pi_1(x)-\frac\beta\gamma\ln\left(\int_{\mathbb{R}} e^{-\gamma  L_{\pi_2}{\bf 0}(x,\pi_1(x),f(x,\pi_1(x),z))}\nu(dz)\right)\\
&=& \pi_1(x)-\frac\beta\gamma\ln\left(\int_{\mathbb{R}} e^{-\gamma  \pi_2(x,\pi_1(x),f(x,\pi_1(x),z))}\nu(dz)\right)\\
&=&\pi_1(x)-\frac\beta\gamma\ln\left(
\int_{\pi_1(x)-x}^\infty e^{-\gamma \pi_2(x,\pi_1(x),x-\pi_1(x)+z)}\nu(dz)+\nu(-\infty,\pi_1(x)-x)\right) .\\
\end{eqnarray*}
Observe that  by  (P1) and the fact that $\pi_k(h_k)\ge 0$ for all $h_k\in H_k$ and $k\in\N,$
 it follows that the sequence $(J_N(x,\pi))_{N\in\mathbb{N}}$ is non-decreasing and
bounded from below by $0$   for every $x\in{\mathbb R}_+$ and $\pi\in\Pi.$
Moreover, for $x\in\mathbb{R}_+,$ $ \pi\in\Pi$ and  $N\in\mathbb{N}$ it holds
\begin{equation}
\label{barb}
J_N(x,\pi)\le x+\bar{b}.
\end{equation}
Indeed, note first that $L_{\pi_N}{\bf 0}(h_N)=\pi_N(h_N)\le x_N+\bar{b}$  for $h_N\in H_N$ with $x_N\ge 0$ and $L_{\pi_N}{\bf 0}(h_N)=0,$
if $x_N<0.$  If $x_{N-1}\ge 0,$ then  making use of (\ref{L_ogr}),  (P3) and (A1) we obtain
\begin{eqnarray*}
 \lefteqn{L_{\pi_{N-1}} (L_{\pi_N} {\bf 0})(h_{N-1})}\\
 &\le&
 \pi_{N-1}(h_{N-1})+\beta \bar{b}-\frac\beta\gamma\ln\left(  \int_{\R} e^{-\gamma f(x_{N-1},\pi_{N-1}(h_{N-1}),z)}\nu(dz)\right)\\
 &\le&
  \pi_{N-1}(h_{N-1})+\beta \bar{b}+\beta \int_{\R} f(x_{N-1},\pi_{N-1}(h_{N-1}),z)\nu(dz)\\
 &\le&
  \pi_{N-1}(h_{N-1})+\beta \bar{b} +\beta(x_{N-1}-\pi_{N-1}(h_{N-1}))+\beta  \int_{0}^\infty z\nu(dz)\\
  &\le&
  \sup_{a\in[0,x_{N-1}]}
 \left( a+\beta(x_{N-1}-a)\right)+\beta \bar{b} +\beta  \mathbb{E}Z^+
=  x_{N-1}+ \bar{b}.
\end{eqnarray*}
If, on the other hand, $x_{N-1}<0,$ then $x_N=x_{N-1}$ and  $L_{\pi_{N-1}} (L_{\pi_N} {\bf 0})(h_{N-1})=0.$
Continuing this procedure and applying (\ref{J_T}), we get the conclusion.
By the above discussion, $ \lim_{N\to\infty} J_N(x,\pi)$ exists for every $x\in \mathbb{R}_+$ and $\pi\in\Pi.$

For an initial level of the risk reserve  $x\in\mathbb{R}_+$ and a policy $\pi\in\Pi,$
we define the risk adjusted discounted cash flow of dividends in the infinite time horizon as follows
\begin{equation}
\label{J}
J(x,\pi):= \lim_{N\to\infty} J_N(x,\pi).
\end{equation}
The aim of the insurance company is to find an optimal value (the so-called value function) of
the risk adjusted discounted cash flow of dividends in the finite and  infinite time horizon, i.e.,
$$ J_N(x):=\sup_{\pi\in\Pi} J_N(x,\pi)\quad\mbox{for  $N\in\N,\ $ and }\quad J(x):=\sup_{\pi\in\Pi} J(x,\pi)$$ and  policies
$\pi_*,\pi^*\in \Pi$ for which
$$J_N(x,\pi_*)=J_N(x)\quad\mbox{for  $N\in\N,\ $ and }\quad J(x,\pi^*)=J(x), \quad \mbox{for all } x\in{\mathbb R}_+.$$

\begin{rmk}\label{r1}
The parameter $\gamma$ represents the risk aversion of the shareholders. The larger $\gamma$, the more risk averse they are.
The limit $\gamma\to 0^+$ leads to the risk neutral case, since
$$ -\frac 1\gamma\ln \Big( \int_{\mathbb{R}_+} e^{-\gamma x}\mu(dx)\Big) \to \int_{\R_+} x\mu(dx),\quad\mbox{for } \gamma\to 0^+. $$
\end{rmk}


\section{The Finite Time Horizon Problem}\label{sec:finite}

In this section, we consider the finite time horizon model. With this end in view we fix the  time horizon, say $N\in\N,$ and  by $V_n$ we denote
 the value function for the problem from period $n$ up to $N,$ where $n=1,...,N,$ i.e.,
$$V_n(h_n)=\sup_{\pi\in\Pi}(L_{\pi_n}\circ \ldots \circ L_{\pi_{N}} ){\bf 0}(h_n),\quad h_n\in H_n.$$
Furthermore,  for $\bar{b}$ defined in (\ref{eq:bbar}), we introduce the set
 \begin{eqnarray*} \mathcal{S}_0 &:=& \{ v:\R\mapsto\R_+ |\ v(x)\le x+\bar{b} \mbox{ for } x\in\R_+, \\
   && \hspace*{1cm}  v(x)=0  \mbox{ for } x<0, \ v \mbox{ is non-decreasing and continuous on }\R_+  \}.
 \end{eqnarray*}
For $v\in\mathcal{S}_0$ we also define the operator $T$ as follows
\begin{eqnarray*}
 \label{T}
  T v(x) &:=&  \sup_{a\in [0,x]} \Big\{a- \frac\beta\gamma \ln \Big( \int_{\R} e^{-\gamma v(f(x,a,z))}\nu(dz)\Big) \Big\}\\ \nonumber
   &=&
   \sup_{a\in [0,x]} \Big\{a- \frac\beta\gamma \ln \Big( \int_{a-x}^\infty e^{-\gamma v(x-a+z)}\nu(dz)+\nu(-\infty,a-x)\Big) \Big\},\ x\in\R_+
\end{eqnarray*}
and
$$ Tv(x)=0, \ x<0.$$
Note that every Borel measurable function $v:\R\mapsto \R_+$ such that $v(x)\le x+\bar{b}$ for $x\in\R_+$ and $v(x)=0$ for $x<0$
can be viewed as a function defined on $H_k,$ with $k\in\N,$ in the sense that
$v(h_k):=v(x_k)$ for every $h_k\in H_k.$ Therefore, with a little abuse of notation,  for any decision rule $\alpha\in\Lambda,$
we shall write
$$ L_\alpha v(x) = \alpha(x) - \frac\beta\gamma \ln \Big( \int_{\alpha(x)-x}^\infty e^{-\gamma v(x-\alpha(x)+z)}\nu(dz)+\nu(-\infty,\alpha(x)-x)\Big) \Big\},\ x\in\R_+$$
and
$$ L_\alpha v(x)=0, \ x<0.$$
We have the following result.

\begin{lem}\label{lem:rewadit}
For any $v\in{\cal S}_0$  it follows that  $Tv\in {\cal S}_0.$
\end{lem}

\begin{pf}
Assume that $x\in\R_+.$ Then,  the continuity  of $Tv$ on $\R_+$ follows from Theorem 2.4.10 in \cite{bra},
since $\mathbb{A}(x)$ is compact, $x\mapsto \mathbb{A}(x)$ is continuous and the mapping
$$(x,a) \mapsto \int_{a-x}^\infty e^{-\gamma v(x-a+z)}\nu(dz)+\nu(-\infty,a-x)$$ is continuous.

We show next that $Tv$ is non-decreasing. Suppose $0\le x_1< x_2$, then we obtain (since $v\ge 0$):
\begin{eqnarray*}
 && \int_{a-x_1}^\infty e^{-\gamma v(x_1-a+z)}\nu(dz)+\nu(-\infty,a-x_1)\\
 &\ge & \int_{a-x_1}^\infty e^{-\gamma v(x_2-a+z)}\nu(dz)+\nu(-\infty,a-x_1)\\
 &=& \int_{a-x_1}^\infty e^{-\gamma v(x_2-a+z)}\nu(dz)+\nu(-\infty,a-x_2)+\int_{a-x_2}^{a-x_1} 1 \nu(dz)\\
 &\ge& \int_{a-x_1}^\infty e^{-\gamma v(x_2-a+z)}\nu(dz)+\nu(-\infty,a-x_2)+\int_{a-x_2}^{a-x_1} e^{-\gamma v(x_2-a+z)} \nu(dz)\\
 &=& \int_{a-x_2}^\infty e^{-\gamma v(x_2-a+z)}\nu(dz)+\nu(-\infty,a-x_2).
\end{eqnarray*}
Taking care of the monotonicity of the logarithm and the minus sign in front, we obtain that
\begin{eqnarray*}
&& a- \frac\beta\gamma \ln \Big( \int_{a-x_1}^\infty e^{-\gamma v(x_1-a+z)}\nu(dz)+\nu(-\infty,a-x_1)\Big)\\
&\le & a- \frac\beta\gamma \ln \Big( \int_{a-x_2}^\infty e^{-\gamma v(x_2-a+z)}\nu(dz)+\nu(-\infty,a-x_2)\Big).
\end{eqnarray*}
The remaining part is a consequence of $[0,x_1]\subset [0,x_2]$.

Finally we prove the upper bound. Setting $u:=x-a$ and making use again of (P1), (P3) and (A1), we conclude
\begin{eqnarray*}
  Tv(x) &=&  \sup_{a\in [0,x]} \Big\{a- \frac\beta\gamma \ln \Big( \int_{a-x}^\infty e^{-\gamma v(x-a+z)}\nu(dz)+\nu(-\infty,a-x)\Big) \Big\} \\
   &=& x+  \sup_{u\in [0,x]} \Big\{-u- \frac\beta\gamma \ln \Big( \int_{-u}^\infty e^{-\gamma v(u+z)}\nu(dz)+\nu(-\infty,-u)\Big) \Big\}\\
   &\le& x+  \sup_{u\in [0,x]} \Big\{-u- \frac\beta\gamma \ln \Big( \int_{-u}^\infty e^{-\gamma (u+z+\bar{b})}\nu(dz)+e^{-\gamma\bar{b}}\nu(-\infty,-u)\Big) \Big\}\\
    &\le& x+ \beta\bar{b}+ \sup_{u\in [0,x]} \Big\{-u- \frac\beta\gamma \ln \Big( \int_{-u}^\infty e^{-\gamma (u+z)}\nu(dz)+\nu(-\infty,-u)\Big) \Big\}\\
   &\le& x+ \beta\bar{b}+ \sup_{u\in [0,x]} \Big\{-u+\beta \int_{-u}^\infty (u+z)\nu(dz) \Big\}\\
      &\le& x+ \beta\bar{b}+ \sup_{u\in [0,x]} \Big\{-u+\beta u +\beta\int_{0}^\infty z\nu(dz) \Big\}=x+\bar{b}.
  \end{eqnarray*}
Clearly,  setting $a:=0$ we also have
\begin{eqnarray*} Tv(x) &\ge& - \frac\beta\gamma \ln \Big( \int_{-x}^\infty e^{-\gamma v(x+z)}\nu(dz)+\nu(-\infty,-x)\Big)\\
 &\ge&  - \frac\beta\gamma \ln \Big( \int_{-x}^\infty \nu(dz)+\nu(-\infty,-x)\Big)= 0.
 \end{eqnarray*}
Hence, the assertion is proved.\hfill$\Box$
\end{pf}

The main result of this section proves  the value iteration for $V_n$ and states that
 the optimal dividend policy is Markov for the model with a finite time horizon.

\begin{thm}\label{theo:finitehorizon}
For every $n=1,\ldots ,N$ we have that $V_{N-n+1}(h_{N-n+1})=J_n(x_{N-n+1})$ and there exists
$\alpha_{N-n+1}^*\in\Lambda$ such that
$J_{n+1}=TJ_n=L_{\alpha_{N-n+1}^*} J_n$, where $J_{0}\equiv {\bf 0}.$ In particular, $J_n\in\mathcal{S}_0$ and $J_n(x)\ge x$ for $x\in \R_+$ and $n\ge 1.$
Moreover, the  policy $\pi_*=(\alpha_1^*,\ldots,\alpha_{N}^*)\in\Pi^M$ is optimal, i.e., $J_N(x)=J_N(x,\pi_*)$
for $x\in\R_+.$
\end{thm}

\begin{pf}
The proof  proceeds by backward induction. Let  $h_N=(x_1,\ldots,x_N)\in H_N$. Then, if $x_N\ge 0$ we obtain
$$ V_{N}(h_N)= \sup_{\pi_N} (L_{\pi_N}{\bf 0})(h_N)= \sup_{a\in [0,x_N]} a = x_N=J_1(x_N)=(TJ_0)(x_N).$$
For $x_N<0$ we put $J_1(x_N)=0.$ Hence, $J_1\in  \mathcal{S}_0.$
Define $\alpha_N^*(x):=x$ for $x\ge0$ and $\alpha_N^*(x):=0$ for $x <0.$ Obviously, $\alpha_N^*\in\Lambda.$
Now suppose that the statement is true  for $k=N,N-1,\ldots,N-n+1,$ ($n\in \N$)  i.e.,
$$V_{N-n+1}(h_{N-n+1}) =J_n(x_{N-n+1}) =(L_{\alpha^*_{N-n+1}}\circ \ldots \circ L_{\alpha^*_{N}} ){\bf 0}(h_{N-n+1}), \quad h_{N-n+1}\in H_{N-n+1}.$$
We prove the result for $k=N-n.$ Fix a history $h_{N-n}\in H_{N-n}$ and assume that $x_{N-n}\ge 0.$ From (\ref{L_mon}) and our assumption we have
\begin{eqnarray} \nonumber\label{max_n}
V_{N-n}(h_{N-n})&=&\sup_{\pi\in\Pi}(L_{\pi_{N-n}}\circ \ldots \circ L_{\pi_{N}} ){\bf 0}(h_{N-n})\\\nonumber
&\le&\sup_{\pi_{N-n}}(L_{\pi_{N-n}}\circ L_{\alpha^*_{N-n+1}}\circ \ldots \circ L_{\alpha^*_{N}} ){\bf 0}(h_{N-n})\\\nonumber
&=&\sup_{\pi_{N-n}}(L_{\pi_{N-n}}V_{N-n+1})(h_{N-n})\\\nonumber
&=&\sup_{a\in[0,x_{N-n}]} \left\{a-\frac\beta\gamma\ln\left(  \int_{\R} e^{-\gamma J_{n}(f(x_{N-n},a,z))}\nu(dz) \right)\right\}\\
&=& (TJ_n)(x_{N-n}) = (L_{\alpha_{N-n}^*}\circ \ldots \circ L_{\alpha^*_{N}}){\bf 0}(x_{N-n})\\\nonumber
&\le& J_{n+1}(x_{N-n}) \le V_{N-n}(h_{N-n}).
\end{eqnarray}
Hence, we have the equality. Since  $\mathbb{A}(x)$ is compact and the set-valued mapping $x\mapsto \mathbb{A}(x)$ is continuous,
the existence of a maximiser $\alpha_{N-n}^*\in \Lambda$ in (\ref{max_n}) follows from, e.g.,  Proposition 2.4.8 in \cite{bra}.
Assume now that $x_{N-n}<0.$ This means that ruin has happened before or at  the epoch $N-n.$  Then, $\alpha^*_{N-n}(x_{N-n})=\ldots=\alpha_N^*(x_N)=0,$
$V_{N-n}(h_{N-n})=\ldots=V_N(h_N)=0$ and $x_{N-n}=\ldots=x_N.$ From Lemma \ref{lem:rewadit}, it follows that $J_{n+1}=TJ_{n}\in {\cal S}_0.$
In order to  conclude, we put $\pi_*=(\alpha_1^*,\ldots, \alpha_N^*).$ Then, $J_N(x)=J_N(x,\pi_*)$.

Now assume that $x\in \R_+$ and observe that $J_1(x)\ge x.$ Putting $a:=x$
we obtain by (P1)
\begin{eqnarray*} J_{n+1}=TJ_n(x) &\ge& x- \frac\beta\gamma \ln
\Big( \int_{0}^\infty e^{-\gamma J_n(z)}\nu(dz)+\nu(-\infty,0)\Big)\\
 &\ge& x - \frac\beta\gamma \ln \Big( \int_{0}^\infty \nu(dz)+\nu(-\infty,0)\Big)\ge x.
 \end{eqnarray*}
 This fact finishes the proof.
\hfill$\Box$
\end{pf}

\begin{rmk}
For $\gamma\to 0^+$ we obtain the value iteration for the risk neutral insurance company
$$ J_{n+1}(x) =  \sup_{a\in [0,x]} \Big\{a+\beta \int_{a-x}^\infty J_{n}(x-a+z)\nu(dz) \Big\}.$$
\end{rmk}


\section{The  Infinite Time Horizon Model}\label{sec:limit}

From  considerations in Section \ref{sec:mod}, it follows that the sequence $(J_N(x))_{N\in\N}$ is also non-decreasing.
Hence, $J_\infty(x):=\lim_{N\to\infty} J_N(x)$ exists and $x\le J_\infty(x)\le x+\bar{b}$ for every $x\in\R_+.$
We arrive at the first result.

\begin{lem}\label{lem:conv}
It holds that $J(x)=J_\infty(x)$ for $x\in\R_+$ and $J$ is non-decreasing.
\end{lem}

\begin{pf} Let $x\in \R_+$ and $\pi\in\Pi.$
Clearly,  we have  $J_N(x)\ge J_{N}(x,\pi).$ Letting $N\to\infty$
 yields that $J_\infty(x) \ge J(x,\pi).$ Hence, taking the supremum over all policies we obtain $J_\infty(x) \ge J(x)$ for all $x\in\R_+.$
On the other hand, for  fixed $N\in\N$ and all $n\ge N$ we get $J_{n}(x,\pi)\ge J_{N}(x,\pi)$. Thus, $J(x,\pi) \ge J_N(x,\pi),$ which
implies that $J(x)\ge J_N(x)$ and, consequently, $J(x)\ge J_\infty(x)$ for all $x\in\R_+.$ Hence, combining both inequalities together we have
that $J(x)=J_\infty(x)$ for $x\in\R_+.$ Since each $J_N$ is non-decreasing, it follows that
$J$ is non-decreasing. \hfill$\Box$
\end{pf}

The second result is a simple observation. For any policy $\pi=\alpha^\infty\in\Pi^S,$ we shall write $J_\alpha(x)$ instead of $J(x,\alpha^\infty)$
and $J_{N,\alpha}(x)$ instead of $J_N(x,\alpha^\infty).$

\begin{lem} \label{lem:fixed_alpha}
Let $\pi=\alpha^\infty\in \Pi^S.$ Then, $J_\alpha =L_\alpha J_\alpha.$
\end{lem}

\begin{pf} From the definition of $J_{N,\alpha}$
it can be easily concluded that
$$ J_{N,\alpha} = L_\alpha J_{N-1,\alpha}=L_\alpha^N {\bf 0},$$
where $L_\alpha^N$ is the $N$-th composition of the operator $L_\alpha$ with itself.
Letting $N\to\infty$ on both sides and making use of the monotone convergence theorem, we get the conclusion.\hfill$\Box$
\end{pf}

The next main result provides a characterisation of the value function in the infinite time horizon model. Let us define
 \begin{eqnarray*} \mathcal{S}&:=& \{ v:\R\mapsto\R_+ |\ x\le v(x)\le x+\bar{b} \mbox{ for } x\in\R_+, \\
   && \hspace*{1cm}  v(x)=0  \mbox{ for } x<0, \ v \mbox{ is non-decreasing and continuous on }\R_+  \}.
 \end{eqnarray*}

\begin{thm} \label{thm:2}
The risk sensitive value function $J$ of the dividend problem is the unique fixed point of $T$ in $\mathcal{S}$, i.e.,
$$ J(x)= T J(x) = \sup_{a\in [0,x]} \Big\{a- \frac\beta\gamma \ln \Big( \int_{a-x}^\infty e^{-\gamma J(x-a+z)}\nu(dz)+\nu(-\infty,a-x)\Big) \Big\},\; x\in\R_+$$
and $J(x)=0=TJ(x)$ for $x<0.$
Moreover, there exists $\alpha^*\in\Lambda$ such that $J=L_{\alpha^*}J.$
\end{thm}

\begin{pf}
We start with defining the set
$$ \mathcal{B} := \{b:\R_+\to\R_+ |\ b(x) \le \bar{b}, b \mbox{ is continuous on } \R_+\}.$$
Let  $v\in {\cal S}$ and  $x\in\R_+.$ Then, $v(x)$ can be written as $v(x)=x+b(x),$ where $b\in{\cal B}.$
Recall that
\begin{eqnarray*}
  Tv(x) &=&  \sup_{a\in [0,x]} \Big\{a- \frac\beta\gamma \ln \Big( \int_{a-x}^\infty e^{-\gamma v(x-a+z)}\nu(dz)+\nu(-\infty,a-x)\Big) \Big\} \\
   &=& x+  \sup_{u\in [0,x]} \Big\{-u- \frac\beta\gamma \ln \Big( \int_{-u}^\infty e^{-\gamma (u+z+b(u+z))}\nu(dz)+\nu(-\infty,-u)\Big) \Big\}.
    \end{eqnarray*}
Defining  the operator $U$ on ${\cal B}$ as follows
\begin{equation}
\label{U}
Ub(x):=\sup_{u\in [0,x]} \Big\{-u- \frac\beta\gamma \ln \Big( \int_{-u}^\infty e^{-\gamma (u+z+b(u+z))}\nu(dz)+\nu(-\infty,-u)\Big) \Big\},
\end{equation}
we obtain that $Tv(x)=x+Ub(x).$
 We claim that $U:{\cal B}\mapsto{\cal B}.$ Indeed,  by (P3) for $x\in\R_+$
   \begin{eqnarray*}
Ub(x)&\le&\sup_{u\in [0,x]} \Big\{-u+\beta\int_{-u}^\infty  (u+z+\bar{b})\nu(dz)\Big\}\\
   &\le&\sup_{u\in [0,x]} \Big\{-u+\beta u +\beta\bar{b}+\beta \mathbb{E} Z^+ \Big\}= \bar{b}.
  \end{eqnarray*}
Moreover, $Ub(x)\ge0$ by taking $u:=0$ in (\ref{U}).

We equip $\mathcal{B}$ with the supremum norm  $\|\cdot\|_\infty.$ Then,
$( \mathcal{B} ,\|\cdot\|_\infty)$ is complete.
We claim that  $U$ defined in (\ref{U}) is a contraction. With this end in view, let $b,c\in \mathcal{B}.$
Since $b\le c+\|b-c\|_\infty,$ we have
\begin{eqnarray*}
&&  Ub(x) -Uc(x) \le  \beta \sup_{u\in [0,x]} \Big\{- \frac1\gamma \ln \Big( \int_{-u}^\infty e^{-\gamma (u+z+b(u+z))}\nu(dz)+\nu(-\infty,-u)\Big) \\
 && \hspace*{3cm}  +\frac1\gamma \ln \Big( \int_{-u}^\infty e^{-\gamma (u+z+c(u+z))}\nu(dz)+\nu(-\infty,-u)\Big) \Big\} \\
 &&\le  \beta \sup_{u\in [0,x]} \Big\{- \frac1\gamma \ln \Big( \int_{-u}^\infty e^{-\gamma (u+z+c(u+z)+\|c-b\|_\infty)}\nu(dz)
 +e^{-\gamma \|c-b\|_\infty}\nu(-\infty,-u)\Big)\\
   && \hspace*{2cm}  +\frac1\gamma \ln \Big( \int_{-u}^\infty e^{-\gamma (u+z+c(u+z))}\nu(dz)+\nu(-\infty,-u)\Big) \Big\} \\
   &=& \beta \|c-b\|_\infty.
\end{eqnarray*}
Exchanging the roles of $b$ and $c$ we get that $\|Ub-Uc\|_\infty \le \beta \|b-c\|_\infty$.

Next we know by Theorem \ref{theo:finitehorizon} that $J_k\in\mathcal{S},$  for $k\in\N,$ and   $J_k=TJ_{k-1}$ for $k\in\N.$
Hence, there exist functions $b_k\in{\cal B}$ for $k\in\N$ such that $J_k(x)= x+b_k(x),$ $x\in\R_+.$
Putting
$\mathrm{id}(x)=x$, we obtain for $x\in\R_+$
$$J_{k}(x)= x+b_{k}(x) = TJ_{k-1}(x) = T(\mathrm{id}+b_{k-1})(x) = x+Ub_{k-1}(x). $$
This implies that $b_k= Ub_{k-1}$ i.e., the bounded part of the value functions $J_k$ can be iterated with the help of the $U$-operator.
On the other hand, by Banach's fixed point theorem the sequence $(b_k)_{k\in\N}$ converges as $k\to\infty$
to a function $b_o\in\mathcal{B},$ which is the unique fixed point of $U.$ Hence, we infer that $J(x)=x+b_o(x)$  for $x\in\R_+$
and $J\in\mathcal{S}$. Therefore,
$$TJ(x)=T(\mathrm{id}+b_o)(x) =x+Ub_o(x)=x+b_o(x)=J(x)$$
for $x\in \R_+$. Since $J(x)=0$ for $x<0$,
we conclude that  $J$ is the unique fixed point of $T$ in $\mathcal{S}$.

The existence of $\alpha^*\in\Lambda$ follows from Proposition 2.4.8 in \cite{bra}.\hfill$\Box$
\end{pf}

\begin{rmk}
The proof of Theorem \ref{thm:2} which essentially hinges on Banach's fixed point theorem shows that the risk sensitive value function can be approximated by sequences of the form
$$ J= \lim_{n\to\infty} T^n J_0$$
with $J_0 \in \mathcal{S}$ and not only by the special sequence $T^n {\bf 0}$ which appears in Lemma \ref{lem:conv}. 
\end{rmk}


\section{Characterising the Value Function $J$ and its Maximiser $\alpha^*$}\label{sec:value}

In what follows we denote by $\alpha^*\in\Lambda$ the largest maximiser of the right-hand side in the following equation
$$ J(x) = \sup_{a\in [0,x]} \Big\{a- \frac\beta\gamma \ln \Big( \int_{a-x}^\infty e^{-\gamma J(x-a+z)}\nu(dz)+\nu(-\infty,a-x)\Big) \Big\}$$
for  $x\in\R_+.$ Since $J$ is continuous it follows from Remark 2.4.9 in \cite{bra} that $\alpha^*(x)$ is upper semicontinuous in $x$.
The next lemma contains  some properties of $J$ and $\alpha^*$.

\begin{lem}\label{lem:propJ}
\begin{itemize}
  \item[a)] For all $x\ge y\ge 0$ it holds that $J(x)-J(y) \ge x-y$.
  \item[b)] For all $x\in\R_+$ it holds that $J\big( x-\alpha^*(x)\big) = J(x)-\alpha^*(x)$ and $\alpha^*\big(x-\alpha^*(x)\big) =0$.
\end{itemize}
\end{lem}

\begin{pf}
\begin{itemize}
  \item[a)] Let $x\ge y\ge 0$. Then by the change of variable $a':=a- x+y$  we obtain that
  \begin{eqnarray*}
    J(x) &=& \sup_{a\in [0,x]} \Big\{a- \frac\beta\gamma \ln \Big( \int_{a-x}^\infty e^{-\gamma J(x-a+z)}\nu(dz)+\nu(-\infty,a-x)\Big) \Big\} \\
     &=& \max\Big\{ \sup_{a\in [0,x-y]} \Big\{a- \frac\beta\gamma \ln \Big( \int_{a-x}^\infty e^{-\gamma J(x-a+z)}\nu(dz)+\nu(-\infty,a-x)\Big) \Big\} , \\
     && \hspace*{0.5cm} \sup_{a\in [x-y,x]} \Big\{a- \frac\beta\gamma \ln \Big( \int_{a-x}^\infty e^{-\gamma J(x-a+z)}\nu(dz)+\nu(-\infty,a-x)\Big) \Big\} \Big\}\\
     &\ge & x-y + \sup_{a'\in [0,y]} \Big\{a'- \frac\beta\gamma \ln \Big( \int_{a'-y}^\infty e^{-\gamma J(y-a'+z)}\nu(dz)+\nu(-\infty,a'-y)\Big) \Big\}\\
     &=& x-y+J(y)
  \end{eqnarray*}
  and the statement follows.
    \item[b)] Let $x\in\R_+$. Then, $x-\alpha^*(x)\ge 0$ and we get by choosing action $a=0$ that
    $$ J\big(x-\alpha^*(x)\big) \ge - \frac\beta\gamma \ln \Big( \int_{\alpha^*(x)-x}^\infty e^{-\gamma J(x-\alpha^*(x)+z)}\nu(dz)+\nu(-\infty,\alpha^*(x)-x)\Big).$$
    On the other hand, by the definition of $\alpha^*$ we obtain
    $$ J(x) = \alpha^*(x)- \frac\beta\gamma \ln \Big( \int_{\alpha^*(x)-x}^\infty e^{-\gamma J(x-\alpha^*(x)+z)}\nu(dz)+\nu(-\infty,\alpha^*(x)-x)\Big).$$
    Thus, we infer
    \begin{eqnarray*}
      J(x) -\alpha^*(x) &=& - \frac\beta\gamma \ln \Big( \int_{\alpha^*(x)-x}^\infty e^{-\gamma J(x-\alpha^*(x)+z)}\nu(dz)+\nu(-\infty,\alpha^*(x)-x)\Big) \\
       &\le &  J\big(x-\alpha^*(x)\big) \le J(x) -\alpha^*(x),
    \end{eqnarray*}
   where the last inequality follows from part a) by setting $y= x-\alpha^*(x)$. Hence,
   we have equality in the last expression and also $\alpha^*\big(x-\alpha^*(x)\big) =0$. This is also the largest maximizer, since
   $\alpha^*(x)$ is the largest maximizer in state $x$. \hfill$\Box$
\end{itemize}
\end{pf}

Next we show that there exists a finite risk reserve level beyond which it is always optimal to pay down to this level.

\begin{lem}\label{lem:xibound}
Let $\xi := \sup\{ x\in\R_+ |\  \alpha^*(x) = 0\}$. Then $\xi <\infty$ and
$$ \alpha^*(x)=x-\xi, \quad \mbox{\rm for all } x\ge \xi.$$
\end{lem}

\begin{pf}
Let $x\in\R_+$ be such that $\alpha^*(x) =0$. Then,  from Section  \ref{sec:mod} we know that $J(x) \le x+\bar{b}.$ Thus,
\begin{eqnarray*}
  J(x) &=&  -\frac\beta\gamma \ln \Big( \int_{-x}^\infty e^{-\gamma J(x+z)}\nu(dz)+\nu(-\infty,-x)\Big) \\
   &\le & -\frac\beta\gamma \ln \Big( \int_{-x}^\infty e^{-\gamma (x+z+\bar{b})}\nu(dz)+e^{-\gamma (x+\bar{b})}\nu(-\infty,-x)\Big)  \\
   &=& \beta x +\beta \bar{b} - \frac\beta\gamma \ln \Big( \int_{-x}^\infty e^{-\gamma z}\nu(dz)+\nu(-\infty,-x)\Big) \\
   &\le &  \beta x +\beta \bar{b} - \frac\beta\gamma \ln \Big( \int_{0}^\infty e^{-\gamma z}\nu(dz)+\nu(-\infty,0)\Big) \\
   &\le &  \beta x +\beta \bar{b}+\beta \int_0^\infty z\nu(dz) = \beta x + \bar{b}.
\end{eqnarray*}
On the other hand,  $J(x) \ge x$. 
Taking into account these two inequalities we get
$$ x\le \frac{\bar{b}}{1-\beta}<\infty,$$
and $\xi$ has to be finite.

Now let $x\ge \xi$. We know from Lemma \ref{lem:propJ}b  that $\alpha^*\big(x-\alpha^*(x)\big) =0$, hence $x-\alpha^*(x)\le\xi$.
Thus a payment of $\alpha^*(x)-(x-\xi)$ is admissible in state $\xi$ and we infer
\begin{eqnarray*}
  J(\xi) &\ge& \alpha^*(x)-(x-\xi) - \frac\beta\gamma \ln \Big( \int_{\alpha^*(x)-x}^\infty e^{-\gamma J(x-\alpha^*(x)+z)}\nu(dz)+\nu(-\infty,\alpha^*(x)-x)\Big)  \\
   &=& J(x)-(x-\xi)  \ge J(\xi).
\end{eqnarray*}
Hence, we have equality and $\alpha^*(x)-(x-\xi)$ is a maximum point in state $\xi$. Since $\alpha^*(\xi)$ is the largest maximum point we obtain
$$0=\alpha^*(\xi) \ge \alpha^*(x)-(x-\xi) \ge 0,$$
which implies that $\alpha^*(x) = x-\xi$.\hfill$\Box$
\end{pf}

Next we will  further characterise $\alpha^*$ on the interval $[0,\xi]$.
It turns out that $\alpha^{*\infty}\in\Pi^S$ is a so-called {\em band policy}.

\begin{defi}\label{def:band}
\begin{itemize}
  \item[a)] A stationary policy $\alpha^\infty$ is called a {\em barrier policy}, if
there exists a number $c\ge 0$ such that
$$\alpha(x) = \left\{ \begin{array}{cl}
0 , & \;\mbox{if}\; x\le c\\
x-c , & \;\mbox{if}\; x>c.
\end{array}\right.$$
\item[b)]  A stationary policy $\alpha^\infty$ is called a {\em (finite) band policy}, if there
exist finitely many numbers $0\le c_0<d_1\le c_1 < d_2\le \ldots < c_m$ such that
$$ \alpha(x) = \left\{ \begin{array}{cl}
0 , & \;\mbox{if}\; x\le c_0\\
x-c_k , & \;\mbox{if}\; c_k<x\le d_{k+1}\\
0 , & \;\mbox{if}\; d_{k+1}< x\le c_{k+1}\\
x-c_m , & \;\mbox{if}\; x > c_m,
\end{array}\right.\quad k\in\N.$$
\item[c)]  A stationary policy $\alpha^\infty$ is called a {\em (countable) band policy}, if there exists a partition of $\mathbb{R}_+$ of the form
$A\cup B=\mathbb{R}_+$ with
$$ f(x) = \left\{ \begin{array}{cl}
 0 & , \mbox{ if } x\in B\\
 x-z, \mbox{ where } z=\sup\{y\; |\; y\in B, 0\le y<x\} & , \mbox{ if } x\in A
 \end{array}\right.$$
\end{itemize}
\end{defi}

Note that a barrier policy is a special finite band policy and a finite band policy is a special countable band policy. In what follows the term 'band-policy' refers to the most general definition in part c).

\begin{thm}
The stationary policy $\alpha^{*\infty}$ is a band policy.
\end{thm}

\begin{pf}
We only have to consider the interval $[0,\xi),$ since  $\alpha^*$ is given
on $[\xi,\infty)$ by Lemma \ref{lem:xibound}.  Let us introduce the function $\Gamma:\R_+\mapsto \R$
$$ \Gamma(x) := -\frac\beta\gamma \ln\Big(\int_{-x}^\infty e^{-\gamma J(x+z)} \nu(dz) + \nu(-\infty,-x) \Big).$$
Next observe that for $0\le y<x\le\xi$  by Lemma \ref{lem:propJ}a, we have
\begin{equation}\label{a}
 J(x) =\sup_{a\in[0,x]}\{ a+\Gamma(x-a)\} \ge  x-y +\sup_{a\in[0,y]}\{a+\Gamma(y-a)\}=x-y+J(y).
\end{equation}
In particular, if $\alpha^*(x)\ge x-y,$ then the action $\alpha^*(x)-x+y\ge 0$ is available in $y.$ Therefore,
from (\ref{a}) it follows that
\begin{eqnarray*}
J(x)&=&\alpha^*(x)+\Gamma(x-\alpha^*(x))\ge x-y +\alpha^*(y)+\Gamma(y-\alpha^*(y))=J(y)+x-y\\
&\ge& x-y+ \alpha^*(x)-x+y +\Gamma(y-(\alpha^*(x)-x+y))=J(x).
\end{eqnarray*}
This implies that all inequalities in the above display become equalities. Since $\alpha^*(y)$ is the largest maximiser in $y$, then
$\alpha^*(y)\ge\alpha^*(x)-x+y$. Assume that  $\alpha^*(y)>\alpha^*(x)-x+y$. Then, for the action $\alpha^*(y)+x-y,$ available in state $x,$
we obtain
$$
J(x)=\alpha^*(x)+\Gamma(x-\alpha^*(x))> x-y +\alpha^*(y)+\Gamma(x-(x-y+\alpha^*(y)))=x-y+J(y)=J(x).
$$
Hence, $\alpha^*(x)=\alpha^*(y)+x-y$.
This fact  can be
used to construct the bands as follows: Let $\alpha^*(x') := \sup_{0\le
x\le \xi}\alpha^*(x)$. The maximal value is attained since $\alpha^*$ is upper
semicontinuous. If $\alpha^*(x')=0$ we are done. Now suppose that
$\alpha^*(x')>0$. Consider the interval $[x'-\alpha^*(x'), x']$. We have $\alpha^*\big(x'-\alpha^*(x')\big)=0$ and
it holds for $x\in [x'-\alpha^*(x'), x']$  that $\alpha^*(x')=\alpha^*(x)+x'-x$. Rewriting this equation as  $\alpha^*(x)=x-\big(x'-\alpha^*(x')\big)$
shows that we have constructed one band of the band policy.  Then we look
for the next highest value on the remaining set $[0,\xi]\setminus
[x'-\alpha^*(x'), x']$. This procedure is carried on until all bands are
constructed. Since every such interval contains at least one rational number and the intervals are disjoint, there
are at most a countable number of them.\hfill$\Box$
\end{pf}


\section{Optimality of $\alpha^*$}\label{sec:optimal}

In this section, we finally show that the stationary policy $\alpha^{*\infty}$  is  optimal  in the infinite time horizon model.

\begin{thm}
The policy $\alpha^{*\infty}\in\Pi^S$  is optimal.
\end{thm}

\begin{pf}
Let  $\alpha\in\Lambda$ and  $v\in\mathcal{S}.$ Then, by
 (\ref{L_ogr})
for some constant $c\in\R_+$  it holds
$L_\alpha (v+c) \le \beta c +L_\alpha v,$
and by induction it follows that $L_\alpha^n (v+c) \le \beta^n c +L_\alpha^n v,$
where $L_\alpha^n$ is the $n$-th composition of the operator $L_\alpha$ with itself. Additionally,  for $\alpha^*$
we have
$$
 L_{\alpha^*} {\bf 0}(x) = \alpha^*(x) \ge (x-\xi)^+ =: p_\xi(x).
$$
Let $x\in\R_+.$ Recalling  that $\mathrm{id}(x)=x$ and making use of  Theorem \ref{thm:2} we infer that
\begin{eqnarray*}
  J(x) &=& L_{\alpha^*}^n J(x) \le L_{\alpha^*}^n (\mathrm{id}+\bar{b})(x) \le L_{\alpha^*}^n (p_\xi +\xi+\bar{b})(x)  \\
   &\le& \beta^n (\xi+\bar{b}) + L_{\alpha^*}^n p_\xi(x) \le \beta^n (\xi+\bar{b}) + L_{\alpha^*}^{n+1}{\bf 0}(x).
\end{eqnarray*}
Letting $n\to\infty$ we obtain $J(x)\le J_{\alpha^*}(x),$ $x\in\R_+.$ However, $J_{\alpha^*}(x)\le J(x).$   Hence, $\alpha^{*\infty}$ is  optimal.\hfill$\Box$
\end{pf}

\begin{thm}
Suppose that the density $g$ is continuously differentiable on the interior of its support.
Then, the value function $J$ is differentiable on $\R_+$ a.e. and $J'\ge1$ a.e.
\end{thm}

\begin{pf}
Recall the structure of the band policy and denote by $I_k=(c_k,d_{k+1})$ the open interval of points, where $\alpha^*(x) = x-c_k$.
From the fixed point equation we obtain for $x\in I_k$
$$ J(x) = x-c_k- \frac\beta\gamma \ln \Big( \int_{-c_k}^\infty e^{-\gamma J(c_k+z)}\nu(dz)+\nu(-\infty,-c_k)\Big) $$
and then $J(x)$ is obviously differentiable with derivative $J'(x)=1$.
Next let $B:= \{x\in\R_+ : \alpha^*(x)=0\}$ and take an interior point $x\in B$. We have
$$  J(x) =  -\frac\beta\gamma \ln \Big( \int_{-x}^\infty e^{-\gamma J(x+z)}g(z)dz+\int_{-\infty}^{-x}g(z) dz\Big)$$
and  using the change of variables $w:=x+z,$ it follows that
$$  J(x) =  -\frac\beta\gamma \ln \Big( \int_{0}^\infty e^{-\gamma J(w)}g(w-x)dw+G(-x)\Big).$$
Hence, we see that due to our assumptions $J'(x)$ exists.
The points where $J$ might not be differentiable are the endpoints of the countable number of intervals $I_k$ and thus countable.

The fact that $J'(x)\ge 1$ follows from Lemma \ref{lem:propJ}a. \hfill$\Box$
\end{pf}


\section{The Policy Improvement Algorithm}\label{sec:PI}

One way to find an optimal dividend policy is to use the Policy Improvement Algorithm,
which however has to be defined in the right way. We also impose in this section the following additional assumption
\begin{itemize}
\item[(A3')] $\nu$ has a density $g$ with respect to the Lebesgue measure which is a.e.\ continuous.
\end{itemize}
Let us set $\xi^* := \frac{\bar{b}}{1-\beta}$
and consider a stationary policy $\alpha^\infty$  such that $\alpha(x) \ge x-\xi^*$ for all $x\ge \xi^*$
and $J_\alpha(x) \ge x$ for $x\in\R_+.$ This is, for example, true for $\alpha(x)=x$. Then,
$J_\alpha(x)=x+\frac{\beta \rho(Z^+)}{1-\beta}$ for $x\in\R_+.$
Now we want to find an improvement of $\alpha$. For this purpose let us again define
\begin{equation}\label{g}
 \Gamma_\alpha(x) := -\frac{\beta}{\gamma} \ln \Big( \int_{-x}^\infty e^{-\gamma J_\alpha(x+z)}\nu(dz) +\nu(-\infty,-x)\Big)
 \end{equation}
and denote by $\delta(x)$ the largest maximiser of
$$ a\mapsto a+\Gamma_\alpha(x-a)$$
on the interval $[0,x]$. Note that such a maximiser exists by Proposition 2.4.8 in \cite{bra}.
Then, it is possible to show that $\delta$ has the following properties.

\begin{thm}\label{thm:pi}
The new decision rule $\delta$ has the following properties:
\begin{itemize}
  \item[a)] $\delta\big( x-\delta(x)\big) = 0$ for all $x\in \R$,
  \item[b)] $\delta(x) \ge x-\xi^*$ for all $x>\xi^*$,
  \item[c)] $x\le J_\alpha(x) \le J_\delta(x) \le x+\bar{b}$ for all $x$.
\end{itemize}
\end{thm}

\begin{pf}
\begin{itemize}
  \item[a)] The statement is true for if $\delta(x)=0$ or $\delta(x)=x$. Suppose now that $0<\delta(x) <x$
  and, on the contrary,  assume that $\delta\big( x-\delta(x)\big) > 0$.
  Thus,   there exists an $a_0\in (0,x-\delta(x)]$ such that
      $$ a_0 + \Gamma_\alpha\big( x-\delta(x)-a_0\big) \ge \Gamma_\alpha(x-\delta(x)).$$
   Since, $\delta$ is the largest maximiser, we have for all $a>\delta(x)$ that
      $$ a+\Gamma_\alpha(x-a) < \delta(x) +\Gamma_\alpha(x-\delta(x)).$$
      Note that $x-\delta(x)-a_0\ge 0$. Combining these two inequalities we obtain:
      \begin{eqnarray*}
        \delta(x) +\Gamma_\alpha(x-\delta(x)) &>& \delta(x)+a_0 + \Gamma_\alpha\big( x-\delta(x)-a_0\big) \\
         &\ge & \delta(x) +\Gamma_\alpha(x-\delta(x)).
      \end{eqnarray*}
      Hence,  $\delta\big( x-\delta(x)\big) = 0$.
    \item[b)] We show first that for $x>\xi^*$ we have $\delta(x)>0$. Consider $a=\alpha(x)$. Here we obtain
    $$\alpha(x) +\Gamma_\alpha(x-\alpha(x)) = J_\alpha(x) \ge x.$$
    For $a=0$ we obtain:
    \begin{eqnarray*}
      \Gamma_\alpha(x) &=& -\frac{\beta}{\gamma} \ln \Big( \int_{-x}^\infty e^{-\gamma J_\alpha(x+z)}\nu(dz) +\nu(-\infty,-x)\Big) \\
       &\le & -\frac{\beta}{\gamma} \ln \Big( \int_{-x}^\infty e^{-\gamma (x+z+\bar{b})}\nu(dz) +e^{-\gamma(x+\bar{b})}\nu(-\infty,-x)\Big) \\
       &=& \beta(x+\bar{b}) -\frac{\beta}{\gamma} \ln \Big( \int_{-x}^\infty e^{-\gamma z}\nu(dz) +\nu(-\infty,-x)\Big) \\
       &\le & \beta(x+\bar{b})+\beta \int_0^\infty z\nu(dz) =\beta x +\bar{b}.
    \end{eqnarray*}
    Hence for $\delta(x)=0$ we necessarily must have that $\beta x+\bar{b} \ge x$ which is the case if and only if $x\le \frac{\bar{b}}{1-\beta} = \xi^*$.
    Thus, for $x>\xi^*$ we must have $\delta(x) >0$. Together with part a) it follows that $\delta(x) \ge x-\xi^*$.
  \item[c)] By definition of $\delta$ we obtain $J_\alpha(x) \le L_\delta J_\alpha(x)$ and by iteration we get
  \begin{eqnarray*}
    J_\alpha(x)  &\le& L^n_\delta J_\alpha(x) \le L^n_\delta(\mathrm{id}+\bar{b})(x) \le
    L^n_\delta\big(p_{\xi^*}+\xi^*+\bar{b}\big)(x) \\
     &\le& \beta^n(\xi^*+\bar{b}) +L^n_\delta p_{\xi^*}(x) \le \beta^n(\xi^*+\bar{b})+ L^{n+1}_\delta{\bf 0}(x),
  \end{eqnarray*}
  where $p_{\xi^*}(x)=(x-\xi^*)^+.$
  Letting $n\to\infty$ the first term on the right-hand side converges to zero and the second term converges to $J_\delta$.
  This implies $x\le J_\alpha(x)\le J_\delta(x) \le x+\bar{b}$ for $x\in\R_+$. \hfill$\Box$
\end{itemize}
\end{pf}

After executing one policy improvement step we obtain a decision rule $\delta$ with a better value $J_\delta$. Repeating this procedure we can further improve the value. In case the improvement step returns the same decision rule, it is optimal. Otherwise we obtain an increasing sequence of value functions which converge against the optimal one.

\begin{thm}
\begin{itemize}
  \item[a)] If $\delta =\alpha$ in the algorithm we have $J_\alpha = J$, i.e., the stationary policy $\alpha^\infty$ is optimal.
  \item[b)] In case the algorithm does not stop, it generates a sequence of decision rules $(\delta_k)$ with $\lim_{k\to\infty}J_{\delta_k} =J$.
\end{itemize}
\end{thm}

\begin{pf}
Before we start with the main part of the proof, the following observation in crucial. When we replace the set $\mathcal{S}$ by
   \begin{eqnarray*} \mathcal{S}' &:=& \{ v:\R\mapsto\R_+ |\ x\le v(x)\le x+\bar{b} \mbox{ for } x\in\R_+, \\
   && \hspace*{1cm}  v(x)=0  \mbox{ for } x<0, \ v \mbox{ is measurable }\R_+  \}.
 \end{eqnarray*}
 Then again $T:\mathcal{S}'\to \mathcal{S}'$. This is true since for $v\in \mathcal{S}'$ we still have
 $$(x,a)\mapsto \int_{\R} e^{-\gamma v(x-a+z)}g(z) dz$$ is continuous. This can be seen as follows:
  Changing variables we get the function
$$(x,a)\mapsto \int_{\R} e^{-\gamma v(w)}g(w-x+a) dw.$$
Assume that $x_n\to x$ and $a_n\to a.$ Hence,
$$\Big|\int_{\R} e^{-\gamma v(w)}g(w-x_n+a_n) dw-\int_{\R} e^{-\gamma v(w)}g(w-x+a) dw\Big|\le$$
$$\int_{\R} |e^{-\gamma v(w)}||g(w-x_n+a_n)-g(w-x+a)| dw\le$$
$$\int_{\R} |g(w-x_n+a_n)-g(w-x+a)| dw\to 0.$$
The last convergence follows from the Scheffe Theorem, since by assumption $g(w-x_n+a_n)\to g(w-x+a)$ for almost all $w\in \R.$
Also in Theorem \ref{thm:2} we can replace $\mathcal{B}$ by
$$ \mathcal{B}' := \{b:\R_+\to\R_+ |\ b(x) \le \bar{b}, b \mbox{ is measurable on } \R_+\}.$$
In total, we obtain that under the additional assumption (A3'), $T$ has a unique fixed point on $\mathcal{S}'$ which is $J$.
\begin{itemize}
  \item[a)] If the algorithm returns $\alpha$ we have $TJ_\alpha =L_\delta J_\alpha =L_\alpha J_\alpha = J_\alpha$ and $J_\alpha\in \mathcal{S}'$. Since $J$ is the unique fixed point of $T$ in $\mathcal{S}'$ we obtain $J=J_\alpha$ and the statement follows.
  \item[b)] 
From Theorem \ref{thm:pi}c we know that $J_{\delta_k}$ is increasing in $k$ and thus the limit $\lim_{k\to\infty}J_{\delta_k} = \bar{J}$ exists and $\bar{J}\in \mathcal{S}'$. Since $J_{\delta_k}\le J$ we have $\bar{J}\le J$. From the definition of $\delta_k$ and Lemma 3 we have $$J_{\delta_{k+1}} = L_{\delta_{k+1}} J_{\delta_{k+1}}\ge L_{\delta_{k+1}} J_{\delta_{k}} = TJ_{\delta_k} \ge L_{\delta_{k}} J_{\delta_{k}}=J_{\delta_k}.$$ Taking the limit $k\to\infty$ we obtain with Theorem A.1.6 in \cite{bra} that $\bar{J}=T\bar{J}$. From our previous discussion it follows that $\bar{J}=J$. \hfill$\Box$
\end{itemize}
 \end{pf}


\section{The Infinite Time Horizon Model:  Case Study}\label{sec:IF}
This section deals with a dividend payout model, in which the increments have the following exponential probability density function
\begin{equation}
\label{distr}
g(x) = \left\{ \begin{array}{ll} \lambda e^{\lambda(x-d)}, & x\le d\\
                                0, & x > d.
                                \end{array}\right.
\end{equation}
Then $G(x)=e^{\lambda(x-d)}$ for $x\le d,$  and $G(x)=1$ for $x>d$ is the cumulative distribution.
Clearly, the mean of $Z$ which has the density in (\ref{distr}) is $d-1/\lambda.$ We should have
$\lambda d>1.$

From Theorem \ref{thm:2} it follows that there exists a function $J\in{\cal S}$ such that
$$ J(x)=\sup_{ a\in [0,x]}\left\{a-\frac\beta\gamma\ln\left(\int_{a-x}^{d}e^{-\gamma J(x-a+z)}g(z)dz+G(a-x)\right) \right\}.$$
Simple re-arrangements and the  substitution $u:=x-a$ give
$$ J(x)=x+\sup_{u\in[0, x]}\left\{-u-\frac\beta\gamma\ln\left(\int_{0}^{d+u}e^{-\gamma J(y)}g(y-u)dy+G(-u)\right) \right\}.$$
Proceeding along similar lines as in \cite{a14} we are able to show that in the risk averse setting the optimal policy is of a barrier type.
With this end in view we set
$$h(u):=-u-\frac\beta\gamma\ln\left(\int_{0}^{d+u}e^{-\gamma J(y)}g(y-u)dy+G(-u)\right).$$
Since $J(x)\le x+\bar{b}$ for every $x\in\R_+,$ then it is easy to infer that $h(u)\to-\infty$ when $u\to\infty.$
Moreover, from the form of function $h$ it follows that it is differentiable on $(0,\infty).$ Therefore,
\begin{eqnarray*}
h'(u)&=&-1+\frac\beta\gamma\left(\frac{\int_{0}^{d+u}e^{-\gamma J(y)}\lambda^2e^{\lambda(y-u-d)}dy-\lambda e^{-\gamma J(d+u)}+\lambda G(-u)}
{\int_{0}^{d+u}e^{-\gamma J(y)}\lambda e^{\lambda(y-u-d)}dy+G(-u)}\right)\\
&=&-1+\frac{\beta\lambda}\gamma\left(1- \frac{ e^{-\gamma J(d+u)}}{\int_{0}^{d+u}e^{-\gamma J(y)}\lambda e^{\lambda(y-u-d)}dy+G(-u)}\right).
\end{eqnarray*}
Suppose first that  $\frac\gamma{\beta\lambda}<1.$
Now, we may have either $h'(0^+)>0$ or $h'(0^+)\le 0.$

$\bullet$ Assume first that $h'(0^+)>0,$ i.e.,
$$-1+\frac{\beta\lambda}\gamma\;\left( \frac{e^{-\gamma J(0)/\beta}-e^{-\gamma J(d)}}{e^{-\gamma J(0)/\beta}}\right)>0.$$
Let $p>0$ be the first point at which $h$ has a local maximum, that is, $h'(p)=0.$ Observe that
$$J(p)=-\frac\beta\gamma\ln\left(\int_{0}^{d+p}e^{-\gamma J(y)}g(y-p)dy+G(-p)\right).$$
Making use of these two facts we find that
$$h'(p)=0=-1+\frac{\beta\lambda}\gamma\;\left( \frac{e^{-\gamma J(p)/\beta}-e^{-\gamma J(d+p)}}{e^{-\gamma J(p)/\beta}}\right),$$
which is equivalent to the equality
$$\ln\left(1-\frac{\gamma}{\beta\lambda}\right)+\gamma J(d+p)=\frac{\gamma J(p)}\beta.$$
Moreover, from Lemma \ref{lem:propJ}a, we know that $J(d+p)-J(p)\ge d.$ Hence, it must hold
 \begin{equation}\label{pp}
J(p)\ge \frac{\frac 1\gamma\ln\left(1-\frac{\gamma}{\beta\lambda}\right)+ d}{1/\beta-1}.
\end{equation}
 On the contrary, assume that there exists $q>p$ at which  $h$ has a global maximum.
 Therefore,  for $x$  lying in the left neighborhood of $q$  we have
 \begin{equation}\label{rq}
 J(x)=-\frac\beta\gamma\ln\left(\int_{0}^{d+x}e^{-\gamma J(y)}g(y-x)dy+G(-x)\right)
 \end{equation}
 and, consequently,
 \begin{equation}\label{rqp}
J'(x)=\frac{\beta\lambda}\gamma\;\left( \frac{e^{-\gamma J(x)/\beta}-e^{-\gamma J(d+x)}}{e^{-\gamma J(x)/\beta}}\right).
\end{equation}
 Obviously, we may take such $x\le q$ for which $x+d\ge q.$ Then, making use of (\ref{rq}) with $x:=q$ and the fact that $q$ is the global
 maximum point  we have
 $$J(x+d)=x+d-q+J(q).$$
Since  $J'(q)=1$ we infer from  (\ref{rqp})  that
$$J'(q)=1=\frac{\beta\lambda}\gamma\;\left( \frac{e^{-\gamma J(q)/\beta}-e^{-\gamma (d+J(q))}}{e^{-\gamma J(q)/\beta}}\right).$$
 Thus,
$$ J(q)=\frac{\frac 1\gamma \ln\left(1-\frac{\gamma}{\beta\lambda}\right)+ d}{1/\beta-1}.$$
However, this equality, (\ref{pp}) and Lemma \ref{lem:propJ}a yield that
$J(p)+q-p\le J(q)\le J(p),$ which leads to $q\le p.$ Hence, $p$ must be the global maximum point of the function $h.$
In this case the optimal policy is of a barrier type:
$$ \alpha^*(x) = \left\{ \begin{array}{ll} 0, & x\le p\\
                                x-p, & x > p.
                                \end{array}\right.$$

$\bullet$ Let us now assume that $h'(0^+)\le 0,$ i.e.,
$$-1+\frac{\beta\lambda}\gamma\;\left( \frac{e^{-\gamma J(0)/\beta}-e^{-\gamma J(d)}}{e^{-\gamma J(0)/\beta}}\right)\le 0.$$
This means that
$$\frac{\gamma J(0)}\beta\ge  \ln\left(1-\frac{\gamma}{\beta\lambda}\right)+ J(d).$$
Making use of Lemma \ref{lem:propJ}a, we obtain the necessary condition for $h'(0^+)\le 0$:
\begin{equation}
\label{war}
J(0)\ge \frac{d+\frac1\gamma\ln\left(1-\frac{\gamma}{\beta\lambda}\right)}{1/\beta-1}.
\end{equation}
Assume that $p$ is the global maximum point of the function $h$. Hence, for $x<p$ (sufficiently close to $p$)
 (\ref{rq}) and (\ref{rqp}) hold true. Clearly, we may consider $x<p$ such that $x+d>p.$ Then,
 $J(x)=x-p+J(p)$ for $x\ge p.$
Combining this equality with (\ref{rqp}) we get
\begin{equation}
\label{J_prim}
J'(x)=\frac{\beta\lambda}\gamma\;\left( \frac{e^{-\gamma J(x)/\beta}-e^{-\gamma(x+d-p+ J(p))}}{e^{-\gamma J(x)/\beta}}\right).
\end{equation}
Letting $x\to p^{-},$ applying that $J'(p)=1$ and (\ref{war}) we infer
$$J(p)=\frac{d+\frac1\gamma\ln\left(1-\frac{\gamma}{\beta\lambda}\right)}{1/\beta-1}\le J(0).$$
However, by Lemma \ref{lem:propJ}a it follows that $J(0)\ge J(p) \ge p+J(0)$
Therefore, the global maximum of the function $h$ must be at $u=0.$ In this case the optimal policy is $\alpha^*(x)=x$ for all $x\in\R_+.$

Consider now the case $\frac\gamma{\beta\lambda}\ge 1.$ Inspecting the derivative of $h$ we see that $h'(u)=-1+\frac{\beta\lambda}\gamma(1- f(u))$ where $f> 0$.
Hence, $h'(u)<0$ for all $u\ge 0$ and its maximum is attained at $u=0.$ Hence, the optimal policy is $\alpha^*(x)=x$ for all $x\in\R_+.$

\section{Influence of the Risk Sensitivity Parameter }\label{sec:RS}

In this section, we discuss the influence of the risk coefficient   $\gamma$ on the optimal policy in the model with the finite time horizon (three stages).
We compute the value function with the help of Theorem  \ref{theo:finitehorizon}. When there is only one payment, we obviously have $J_1(x)=x$
independent of $\gamma$. Now consider $J_2$. We obtain by the transformation $u:=x-a$ for $x\in\R_+$ and by plugging in the density $g$
that
\begin{eqnarray*}
  J_2(x) &=&  \sup_{a\in [0,x]} \left\{a- \frac\beta\gamma \ln \Big( \int_{a-x}^\infty e^{-\gamma (x-a+z)}\nu(dz)+\nu(-\infty,a-x)\Big) \right\}\\
   &=& x+\sup_{u\in [0,x]} \left\{-u- \frac\beta\gamma \ln \Big( e^{-\gamma u}\int_{-u}^\infty e^{-\gamma z} g(z)dz
   + \int_{-\infty}^{-u} g(z) dz \Big)\right\}.
\end{eqnarray*}
With a little abuse of notation define the function $h,$ which has to be maximised
$$ h(u) := -u- \frac\beta\gamma \ln \Big( e^{-\gamma u}\int_{-u}^\infty e^{-\gamma z} g(z)dz+
\int_{-\infty}^{-u} g(z) dz \Big) .$$
In order to look for the maximum we differentiate this function and obtain
$$ h'(u) = -1 +\beta \frac{e^{-\gamma u} \int_{-u}^\infty e^{-\gamma z} g(z)dz}{e^{-\gamma u}
\int_{-u}^\infty e^{-\gamma z} g(z)dz +\int_{-\infty}^{-u} g(z)dz}.$$
Since $\beta<1$ and  the density is non-negative, it is easy to see that $h'(u) < 0$ for all $u\ge 0,$
which means that $h$ is decreasing and the maximum is attained at $u=0$.
Being aware of the transformation we obtain $\alpha_2^*(x)=x$, i.e., the optimal decision rule is to pay out everything
at the beginning of a planning horizon of length two, independent of $\gamma$. Hence, we conclude that
\begin{equation}\label{J_2}
 J_2(x) = x -\frac\beta\gamma \ln \Big( \int_{0}^\infty e^{-\gamma z} g(z)dz+ \int_{-\infty}^{0} g(z) dz \Big) = x+\beta \rho(Z^+).
 \end{equation}
In particular, in the risk neutral case we get $J_2(x) = x+\beta \mathbb{E}Z^+$.

Next we consider $J_3$. Making use of (\ref{J_2}) and of  the same transformation as before we get
\begin{eqnarray*}
  J_3(x) &=&  \sup_{a\in [0,x]} \Big\{a- \frac\beta\gamma \ln \Big( \int_{a-x}^\infty e^{-\gamma (x-a+z+\beta \rho(Z^+))}\nu(dz)
  +\nu(-\infty,a-x)\Big) \Big\}\\
   &=& x+\sup_{u\in [0,x]} \Big\{-u- \frac\beta\gamma \ln
   \Big( e^{-\gamma (u+\beta \rho(Z^+))}\int_{-u}^\infty e^{-\gamma z} g(z)dz+ \int_{-\infty}^{-u} g(z) dz \Big)\Big\}.
\end{eqnarray*}
We define, again abusing the notation, the function $h$ as follows
$$  h(u) := -u- \frac\beta\gamma \ln \Big( e^{-\gamma (u+\beta \rho(Z^+))}\int_{-u}^\infty e^{-\gamma z} g(z)dz
+ \int_{-\infty}^{-u} g(z) dz \Big).$$
Differentiating $h$ yields
$$ h'(u) = -1 +\beta \Big( 1- \frac{ \int_{-\infty}^{-u}g(z) dz +\frac1\gamma
\big( 1- e^{-\gamma \beta \rho(Z^+)}\big) g(-u)}{e^{-\gamma (u+\beta \rho(Z^+))} \int_{-u}^\infty e^{-\gamma z} g(z)dz
+\int_{-\infty}^{-u} g(z)dz}\Big).$$
In case of the risk neutral setting  ($\gamma\to 0^+$) the expression is given by
$$ h'(u) = -1 +\beta \int_{-u}^\infty  g(z)dz+\beta^2 \mathbb{E}Z^+ g(-u).$$
Here it is easy to see by inspection of $h''$ that $h'$ is decreasing, if the density $g$ is increasing and log-concave on
$(-\infty,0)$ and $\frac{g(0)}{g'(0)}\le \beta \mathbb{E}Z^+$. Now if
$h'$ is decreasing we can either have $h'(0)\le0$ in which case $h'(u)\le 0$ for all $u$
and the maximum point is again $u=0$ or $h'(0)>0,$ in which case $h$ is first increasing on an interval $[0,q)$
and then decreasing on $(q,\infty)$. Hence, $q$ is the maximum point of $h$ and the optimal dividend payout is a barrier with size $q$.

\begin{figure}
    \centering
      \includegraphics[trim = 6mm 30mm 10mm 10mm, clip, height=8cm]{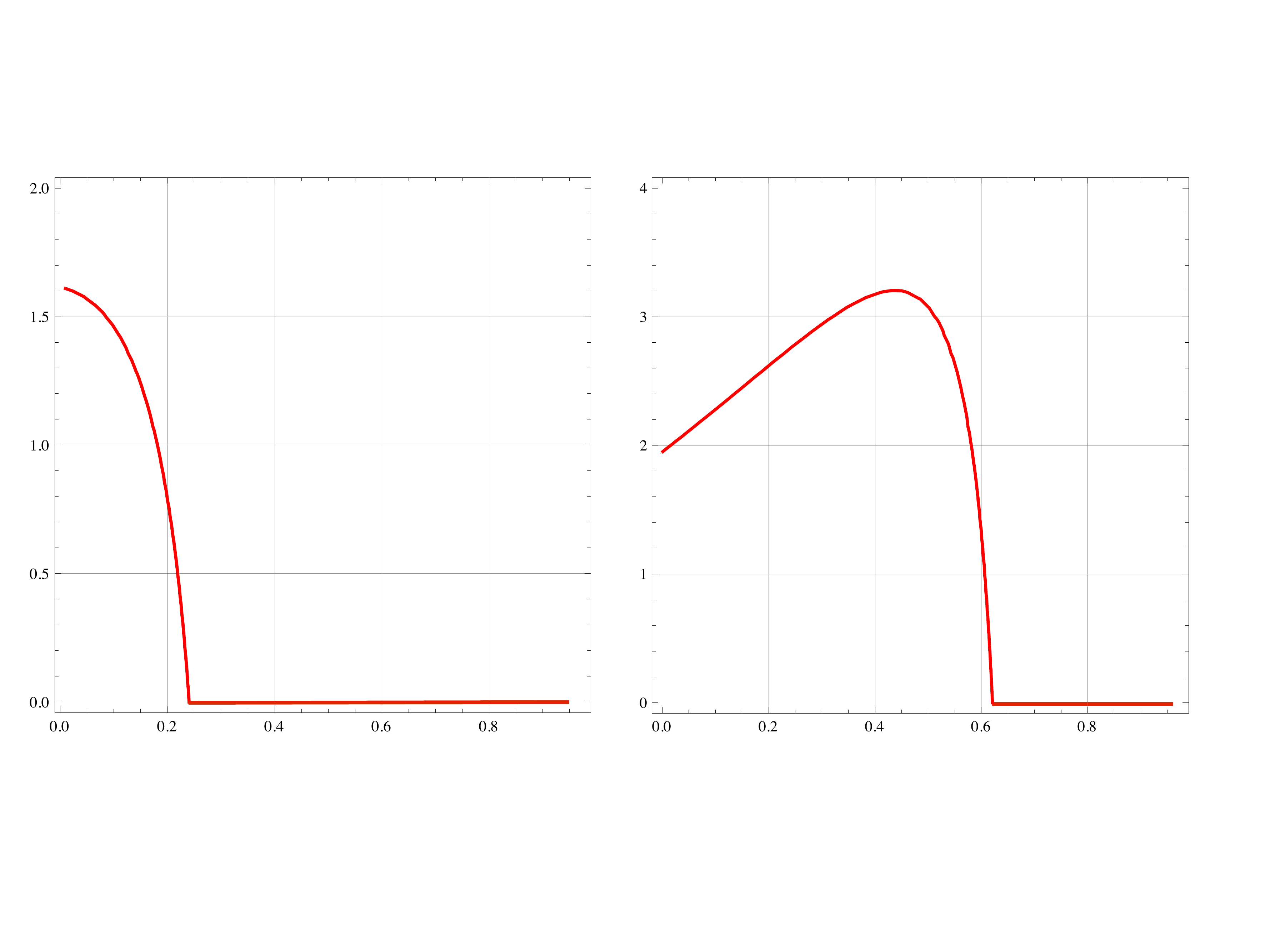}
   \caption{The barrier  as a function of $\gamma$. The left-hand side with $\mu=1.2$. The right-hand side with $\mu=2$.}
    \label{fig:h}
\end{figure}

\begin{figure}
    \centering
         \includegraphics[trim = 6mm 30mm 10mm 10mm, clip, height=8cm]{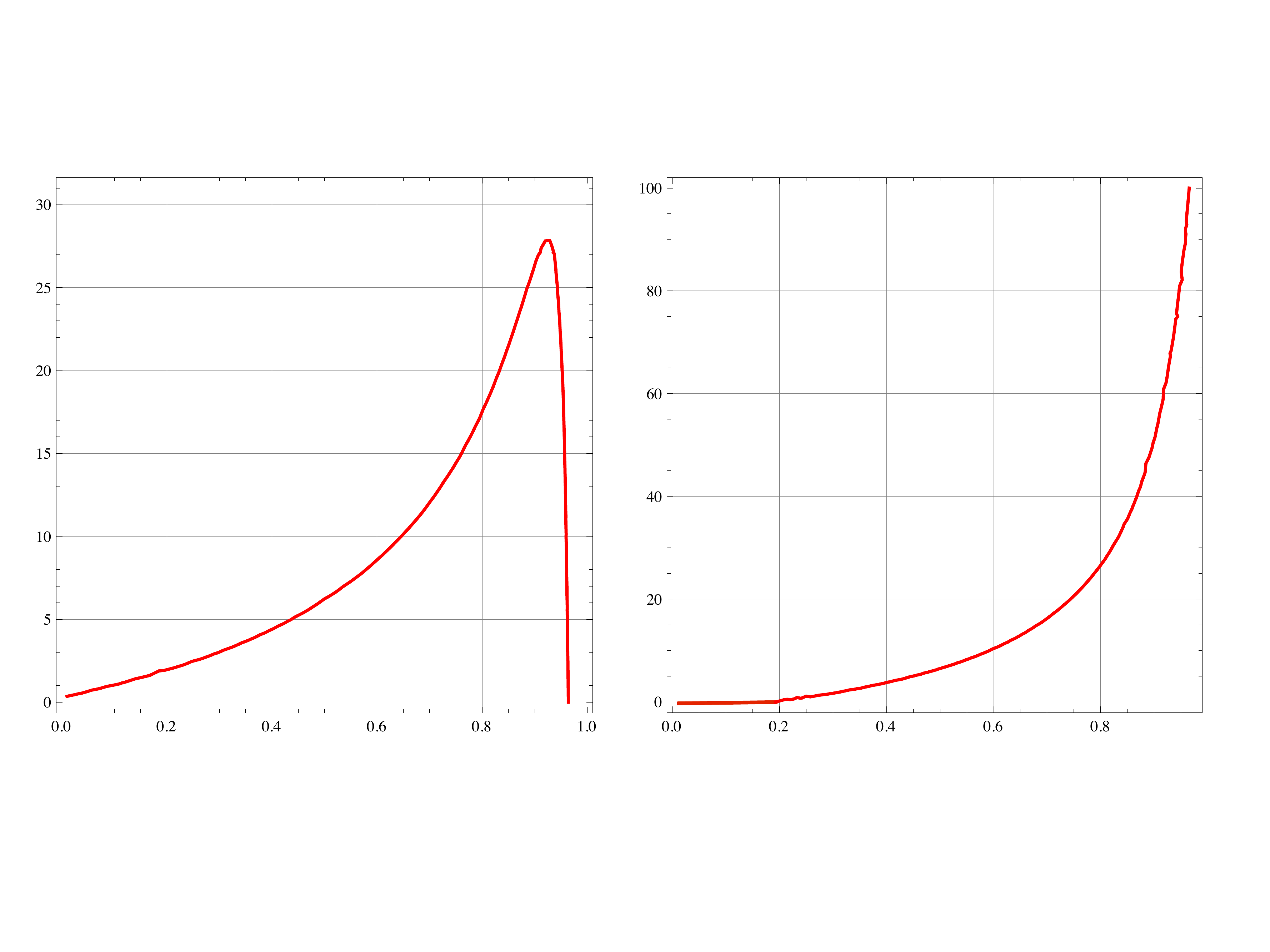}
   \caption{The barrier  as a function of $\gamma$. The left-hand side with $\mu=5$. The right-hand side with $\mu=8$.}
    \label{fig:h2}
\end{figure}

\begin{ex}
Since the risk sensitive case is not so easy to discuss in general, we consider a specific example for the density,
namely the so-called {\em double-exponential} with mean $\mu$, i.e.,
$$ g(x) = \left\{ \begin{array}{ll} \frac12 e^{-(\mu-x)}, & x\le \mu\\
                                \frac12 e^{-(x-\mu)}, & x > \mu
                                \end{array}\right.$$
We have set $\beta=0.99$ in all calculations. In Figures  \ref{fig:h}  and \ref{fig:h2}
we have plotted the barrier as a function of $\gamma$ for different $\mu$.  For $\gamma\to 0^+$ we obtain the risk neutral situation.
The behaviour of this barrier is intriguing. It is very sensitive to the chosen parameter $\mu$, which is the mean of $Z$.
It is worthy to notice that the variance and further central moments are constant and  independent of $\mu.$ Therefore,
we shall discuss the evolution of the curve when the expectation $\mu$ of $Z$ is increasing.
For small values of $\mu$ we can see that the barrier is decreasing,
when $\gamma$ is increasing, i.e., more risk averse shareholders prefer earlier payments.
This may be due to the fear of an early ruin.
However, if the expectation $\mu$ is larger and the company has a good probability to survive for some time period,
 the barrier is first increasing,
i.e., shareholders prefer later payments, which are then rather regular. But surprisingly this is only true up to a certain level of $\gamma$.
Beyond that level, the barrier decreases rapidly  until it gets zero.
This means that very risk averse shareholders prefer to have their money at once. It seems that both payment policies,
where either a very high barrier is set in order to produce a regular dividend stream or the money is paid out at once,
which has also a low variability are reasonable for risk sensitive shareholders. Obviously from an economic point of view
the first policy is more meaningful. Very risk averse shareholders seem to be bad for a company.
\end{ex}
\begin{figure}
    \centering
     \includegraphics[trim = 6mm 30mm 10mm 10mm, clip, height=8cm]{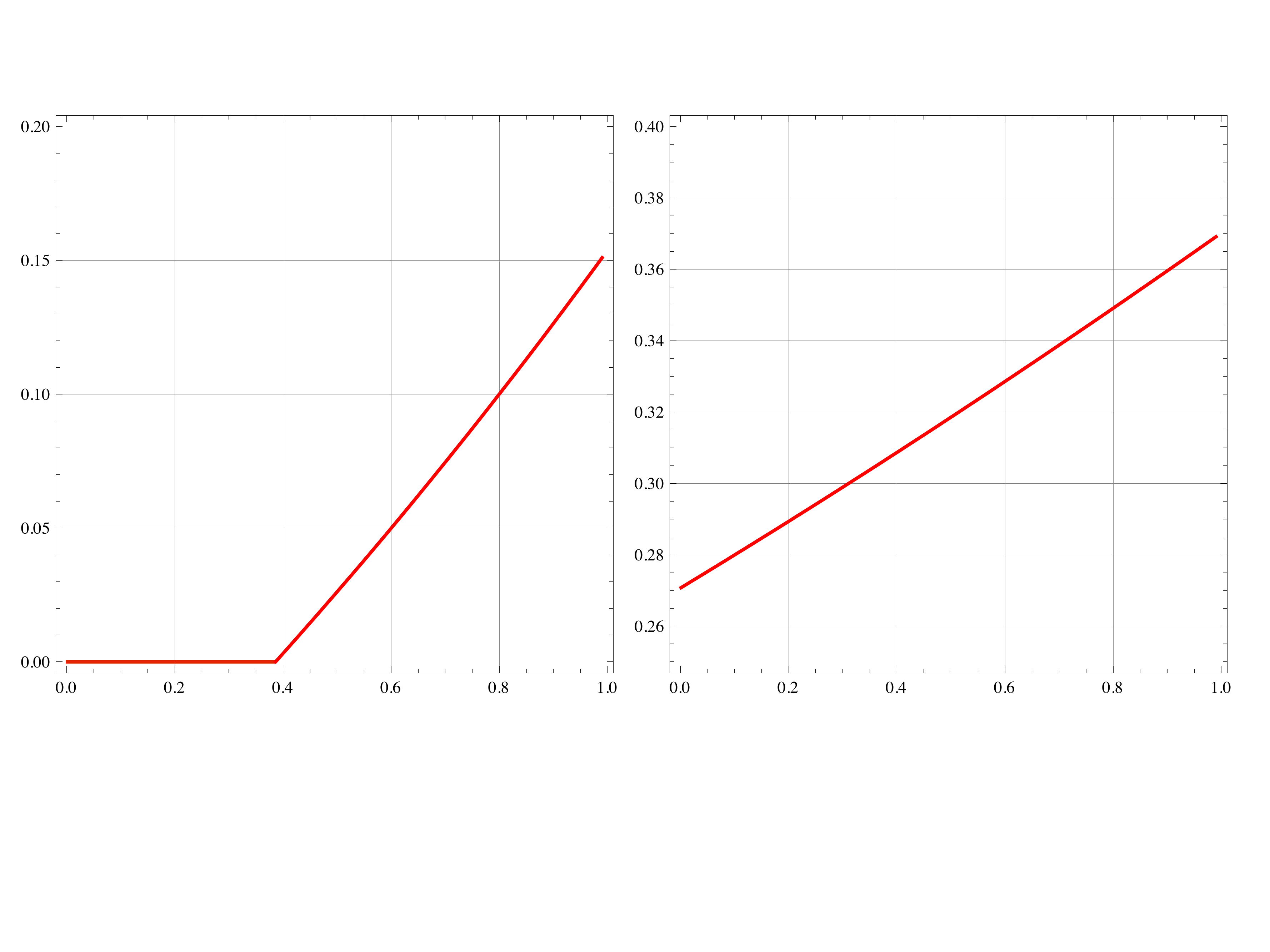}
   \caption{The barrier  as a function of $\gamma,$ if $\lambda=6.$ The left-hand side with $d=1.1.$ The right-hand side
   with $d=0.5$.}
    \label{fig:hp1}
\end{figure}

\begin{figure}
    \centering
     \includegraphics[trim = 6mm 30mm 10mm 10mm, clip, height=8cm]{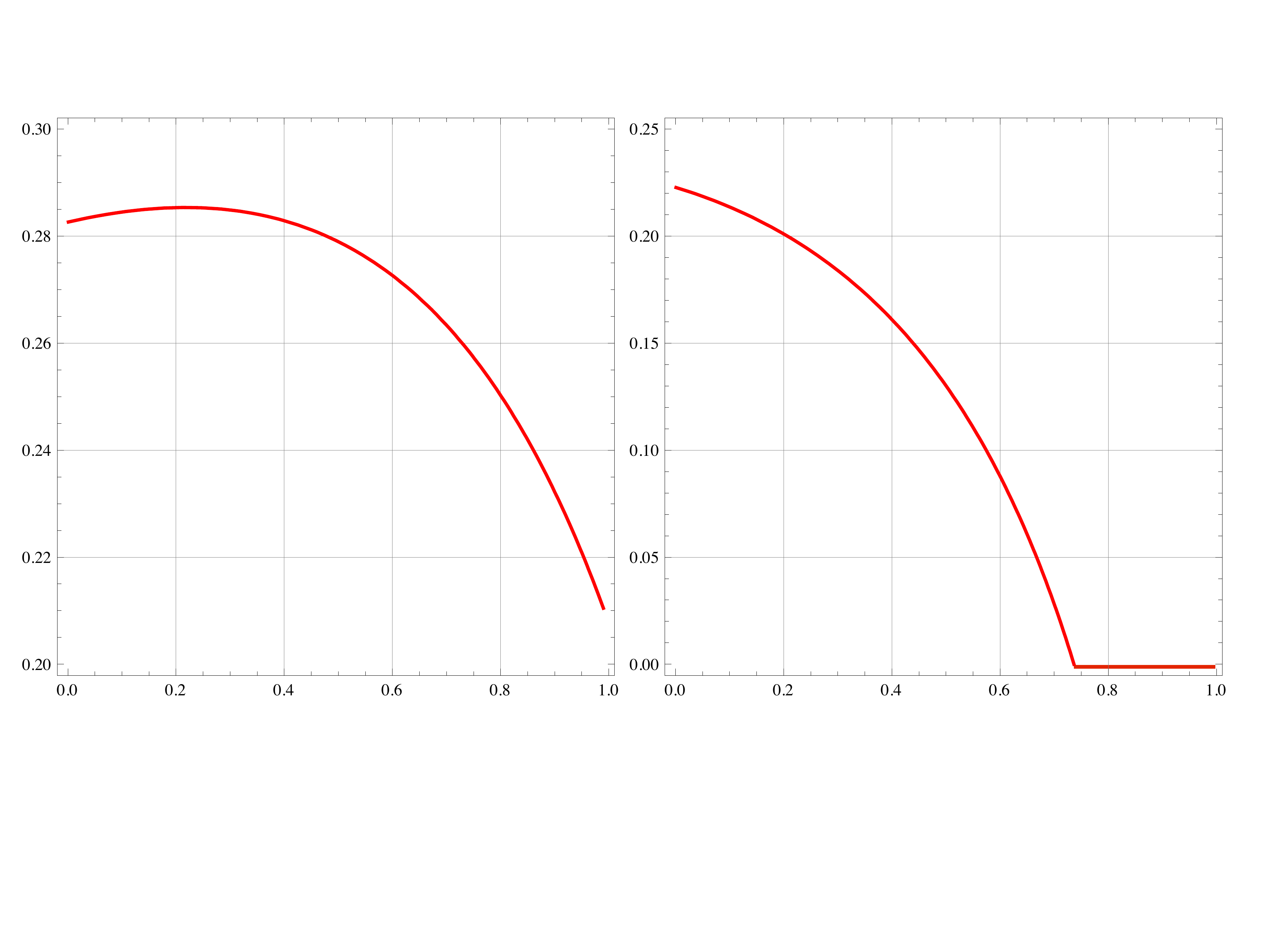}
   \caption{The barrier  as a function of $\gamma,$ if $d=0.5.$ The left-hand side with $\lambda=4.5.$ The right-hand side
   with $\lambda=4.2$.}
    \label{fig:hp2}
\end{figure}

\begin{ex} Let us now consider the distribution defined in (\ref{distr}).  This distribution has the mean equal to $d-1/\lambda$ and
the variance equal to $1/\lambda^2.$
We can see that both the first and second moments play a crucial role in determining the barrier.
In Figure \ref{fig:hp1} the variances of $Z$ are the same, but the means  are different. It can be seen that
the shareholders in case of larger expectation of $Z$
 are willing to get  payments at once. If they are more risk averse than the barrier starts increasing. If the mean of $Z$ is smaller (the second picture in
 Figure \ref{fig:hp1}),
 then the barrier increases at once
 together with the values of risk coefficient. Hence, if the shareholders expect that the risk reverse is stable,
 in the sense that the company will not be ruined so fast,
 they wish to have payments at once.  Otherwise, they prefer to wait
 until the risk reserve attains some critical value. However, the more risk averse shareholders wish to wait longer for their dividends.
This behaviour is in contrast to the case, when the mean of $Z$  is rather small, but the variance of $Z$ is larger.
Figure \ref{fig:hp2} shows that  the barrier decreases,
 either at once or at a certain point, when the decision maker becomes more risk averse.
 This means that the risk neutral shareholders or not too much risk averse shareholders
 prefer to wait for the payments until some critical point.
  If, on the other hand,  they are very risk averse, then they wish
 to have their dividends at once.
\end{ex}

\noindent
{\bf Acknowledgement.} We thank both reviewers for careful reading of the manuscript and their comments that improved the presentation of the paper.

\vspace*{1cm}
\noindent{\bf\large References}

\bibliographystyle{abbrv}

\end{document}